\documentclass[letter,twocolumns,secthm,amsthm,seceqn]{autart}
\usepackage{amsmath}
\usepackage{enumitem}
\usepackage{amssymb}
\usepackage{amsfonts}
\usepackage{cancel}
\usepackage{graphics}
\usepackage{graphicx}
\usepackage{float}
\usepackage{caption}
\usepackage{subcaption}
\usepackage{multicol}
\usepackage{theoremref}
\usepackage{bm}
\usepackage[title]{appendix}
\usepackage{epstopdf}

\begin{document}
	
	\begin{frontmatter}
		\title{Lie Bracket Approximation-Based Extremum Seeking with Vanishing Input Oscillations\thanksref{footnoteinfo}}
		\thanks[footnoteinfo]{Corresponding author M.~Abdelgalil.}
		\author[ADCL]{Mahmoud Abdelgalil}\ead{maabdelg@uci.edu},
		\author[ADCL]{Haithem Taha}\ead{hetaha@uci.edu}
		\address[ADCL]{Department of Mechanical and Aerospace Engineering\\ University of California, Irvine, CA, 92617}

		\begin{abstract}
			In recent years, an approach to extremum seeking control made it possible to design control vector fields that lead to asymptotic stability of the minimum point provided that the minimum value of the function is known a priori. In this work we aim to relax that assumption. We propose an extremum seeking control law that converges to the minimum point with vanishing control oscillations, without access to the minimum value of the cost function. We provide a numerical example to support our results.
		\end{abstract}
		
	\end{frontmatter}

	\section{Introduction}
	Extremum seeking control is an adaptive control technique that drives the steady-state response of a dynamical system to a neighborhood of the minimum point of a cost function in the absence of direct access to gradient information. For more details, the reader is referred to \cite{ariyur2003real,oliveira2017,suttneroutput,tan2010extremum} and the references therein. In this paper, we focus on extremum seeking systems that exploit high amplitude, high frequency, sinusoidal signals. This type of signal is prominently used in motion planning of nonholonomic systems \cite{hassan2020,liu1997b,murray1993}, and techniques from averaging theory are typically applied for analysis \cite{liu1997a,maggia2020}. The first connection to the motion planning framework appears in the reference \cite{Durr2013}. Thenceforth, several authors have contributed to this line of work, e.g. \cite{durr2015singularly,durr2017extremum,grushkovskaya2020extremum,Grushkovskaya2017,labar2018gradient,scheinker2014extremum,scheinker2014non,sutt2019,suttner2020extremum,Suttner2017}.
	
	Traditional extremum seeking \cite{Durr2013,KRSTIC2000595} suffers from persistent oscillations of the steady state response around the minimum point. A solution was proposed in \cite{scheinker2014non}, where the authors extended the averaging techniques in \cite{liu1997a} to nonsmooth systems, which enabled analysis for a set of nonsmooth control vector fields with useful properties such as vanishing at minimum points. A different set of nonsmooth control functions was proposed in \cite{Suttner2017}, which allowed asymptotic stability of the minimum point. Later, it was shown in \cite{Grushkovskaya2017} that both sets of control functions proposed in \cite{scheinker2014non,Suttner2017} belong to a unifying class of generating vector fields. Nevertheless, one of the main assumptions in all these efforts \cite{grushkovskaya2020extremum,Grushkovskaya2017,scheinker2014non,Suttner2017} to guarantee asymptotic convergence to the minimum is that the function value at the minimum point is known a priori. This was pointed out explicitly in several locations, for instance in \cite{grushkovskaya2020extremum,Grushkovskaya2017}.
	
	The contribution of this paper is to provide an extension of the results highlighted so far to the case when the minimum value of the function is unknown. Specifically, we prove asymptotic convergence to the minimum point with bounded amplitude and frequency, for all initial conditions in a subset of the \textit{epigraph} of the cost function.

	\section{Main Theorem}
    \textbf{Notations}: We use bold characters to distinguish vectors and vector valued maps from scalars. Let $D\subset\mathbb{R}^n$ be a subset. We denote the set of $k$-times continuously differentiable real-valued functions on $D$ by $C^k(D;\mathbb{R})$. We denote by $\bm{e}_j$ the $j^{th}$ canonical unit vector in $\mathbb{R}^{n}$. The set of vector fields with regularity $\nu \in \mathbb{N}$ on $\mathbb{R}^n$ is denoted by $\Gamma^\nu(\mathbb{R}^n)$. The Lie derivative of a function $g\in C^1(\mathbb{R}^n;\mathbb{R})$ along a vector field $\bm{f}\in\Gamma^\nu(\mathbb{R}^n)$ is written as $L_{\bm{f}}\, g (\bm{x})$. The Lie bracket between two vector fields $\bm{f}_1,\bm{f}_2 \in\Gamma^1(\mathbb{R}^n) $ is computed as $[\bm{f}_1,\bm{f}_2] (\bm{x})= \textbf{J}_{\bm{f}_2}(\bm{x}) \bm{f}_1 (\bm{x}) - \textbf{J}_{\bm{f}_1}(\bm{x}) \bm{f}_2 (\bm{x}),$ where $\textbf{J}_{\bm{f}}(\bm{x})$ is the standard Jacobian of $\bm{f}$ in the $\bm{x}$-coordinates.
    
	Let $D\subset\mathbb{R}^n$ be a bounded subset with a nonempty interior. Suppose that $J:D\rightarrow\mathbb{R}$ is a cost function that has the following properties:
	\begin{assum}\thlabel{A1}
	Assume that $J\in C^3(D;\mathbb{R})$ and that there exists a unique point $\bm{x}^*\in D$, such that $\tilde{J}(\bm{x}) = J(\bm{x})-J^*>0~ \forall \bm{x}\neq \bm{x}^*$, where $J^* = J(\bm{x}^*)$ and
	\begin{equation*}
	    \begin{aligned}
	        \kappa_1 \tilde{J}(\bm{x})^{2-\frac1m}\leq \lVert\nabla J(\bm{x})\lVert^2\leq \kappa_2 \tilde{J}(\bm{x})^{2-\frac1m}\\
	        \gamma_1 \tilde{J}(\bm{x})^{1-\frac1m}\leq \lVert\nabla^2J(\bm{x})\lVert\leq \gamma_2 \tilde{J}(\bm{x})^{1-\frac1m}
	    \end{aligned}
	\end{equation*}
	where $\kappa_i,\gamma_i>0$ and $m\geq1$.
	\end{assum}

    Next, define the \textit{epigraph} and \textit{strict epigraph} of $J$ 
	\begin{equation*}
    	\begin{aligned}
    	\text{epi}(J)&= \big\{(\bm{x},z)\in D\times \mathbb{R}\big|J(\bm{x})-z\leq0\big\}\\
    	\text{epi}_S(J)&= \big\{(\bm{x},z)\in D\times \mathbb{R}\big|J(\bm{x})-z<0\big\},
    	\end{aligned}
	\end{equation*}
	Let $\bm{\theta}=(\bm{x},z)\in\text{epi}_S(J)$, and define the functions $g_i:\text{epi}_S(J)\rightarrow\mathbb{R},~i\in\{1,2,3\}$
	\begin{equation}\label{g_fcns}
	\begin{aligned}
	g_1(\bm{\theta}) &= \tilde{J}(\bm{x})- J_0,	\quad g_2(\bm{\theta}) = z  - J^* - z_0 \\
	g_3(\bm{\theta}) &= \frac{\tanh(\tilde{J}(\bm{x})^{2-\frac1m})}{z-J(\bm{x})} - y_0
	\end{aligned}
	\end{equation}
	where $J_0,y_0,z_0$ are positive constants. Let $\epsilon> 0$ and define the domains
	\begin{equation*}
	\begin{aligned}
	\Delta_0 &= \big\{ \bm{\theta}\in \text{epi}_S(J) \big| g_i(\bm{\theta})\leq 0, ~\forall i \in \{1,2,3\} \big\}\\
	\Delta_\epsilon &= \big\{ \bm{\theta}\in \text{epi}_S(J) \big| g_i(\bm{\theta})\leq\epsilon, ~\forall i \in \{1,2,3\} \big\}
	\end{aligned}
	\end{equation*}
	 Let $\Lambda$ denote the set of all ordered pairs $(j,s)$, where $ j\in\{1,...,n\}, s\in\{1,2\}$. Then, consider the dynamical system
	\begin{equation}\label{es_sys}
	\begin{aligned}
	\dot{\bm{\theta}}&= \bm{f}_0(\bm{\theta}) + \sum_{\lambda\in\Lambda} \bm{f}_{\lambda}(\bm{\theta}) u_{\lambda}(t)\\
	\bm{f}_0(\bm{\theta}) &= -(z - J(\bm{x}))~\bm{e}_{n+1}\\
	\bm{f}_{j,s}(\bm{\theta}) &= F_s(z-J(\bm{x}))~\bm{e}_j
	\end{aligned}
	\end{equation}
	and the functions $u_{j,1},u_{j,2}$ (the dithers) are given by
	\begin{equation}
	\begin{aligned}
	u_{j,1}(t) &= 2\sqrt{\pi \omega_j \omega}\sin(2\pi\omega_j\omega t) \\
	u_{j,2}(t) &= 2\sqrt{\pi \omega_j \omega}\cos(2\pi\omega_j\omega t)
	\end{aligned}\label{eq:dithers}
	\end{equation}
	where $\omega\in(0,\infty),~\omega_j\in\mathbb{N},~\forall j \in \{1,2,...,n\}$, and the functions $F_{1}(\cdot), F_{2}(\cdot)$ are given by:
	\begin{equation}
	    \begin{aligned}
	        F_{1}(y) &= \sqrt{y}\sin{(\log (y))}\\
	        F_{2}(y) &= \sqrt{y}\cos{(\log( y))}
	    \end{aligned}
	    \label{eq:vecs}
	\end{equation}
	Note that other choices for $F_s(\cdot)$ are possible \cite{Grushkovskaya2017}. Also note that the second equation in system (\ref{es_sys}) is similar to the approach in \cite{sutt2019}.
	\begin{thm}\thlabel{main_thm}
		\normalfont Suppose that the function $J$ satisfies \thref{A1}, and consider the system defined by (\ref{es_sys}), (\ref{eq:dithers}) and (\ref{eq:vecs}). Fix a choice for the collection of frequencies $~\omega_j\in\mathbb{N},~\forall j \in \{1,2,...,n\}$ such that $\forall i\neq j$, $ \omega_i\neq\omega_j$. Then, $\exists\omega^*\in(0,\infty)$ such that $\forall\omega\in(\omega^*,\infty), ~ \forall \bm{\theta}(0)\in \Delta_0$ we have:
		\begin{enumerate}
			\item $\bm{\theta}(t) \in \Delta_\epsilon, \forall t\in[0,\infty)$, 
			\item $\bm{\theta}(t)\rightarrow (\bm{x}^*,J^*) \text{ as }~ t\rightarrow \infty$.
		\end{enumerate}
	\end{thm}
	\begin{pf}
	    The proof is in appendix B. Note that we outline a procedure to estimate a sufficiently high frequency $\omega^*$ in \thref{lem2,lem3}. 
	\end{pf}
	\begin{rem}
		Since $J(x(0))$ is available via measurement, it is always possible to place $z(0)$, which is an internal state of the controller, such that the initial condition strictly lies in $\Delta_0$. We emphasize that this does not require additional information other than online measurement of the function value.
	\end{rem}

	\section{Numerical Simulations}
	\begin{exmp}
	Let $\bm{x}\in\mathbb{R}^2$, and consider the dynamical system 
	\begin{equation}
	    \begin{aligned}
	        \dot{\bm{x}} &= \textbf{A}(t)(\bm{x}-\bm{x}^*) + \textbf{B}\bm{u} \\
	        J(\bm{x},t) &= \frac32\lVert\bm{x}-\bm{x}^*\lVert^2+5(1+\exp(-t))
	    \end{aligned}
	\end{equation}
	where $\bm{u}\in\mathbb{R}^2$ is the control input, $\bm{x}^*=(1,-1)$, and
	$$\textbf{A} (t)=\begin{bmatrix}
	    \cos(t)^2 & \sin(t)^2 \\
	    -\sin(t)^2 & \cos(t)^2
	\end{bmatrix}, \quad 
	\textbf{B}=\begin{bmatrix}
	    1 & 1 \\
	    -1 & 1
	\end{bmatrix}$$
		    \begin{figure}[H]
        	\begin{subfigure}[h]{.49\columnwidth}
                \centering
        	    \includegraphics[width=\textwidth,height=\textwidth]{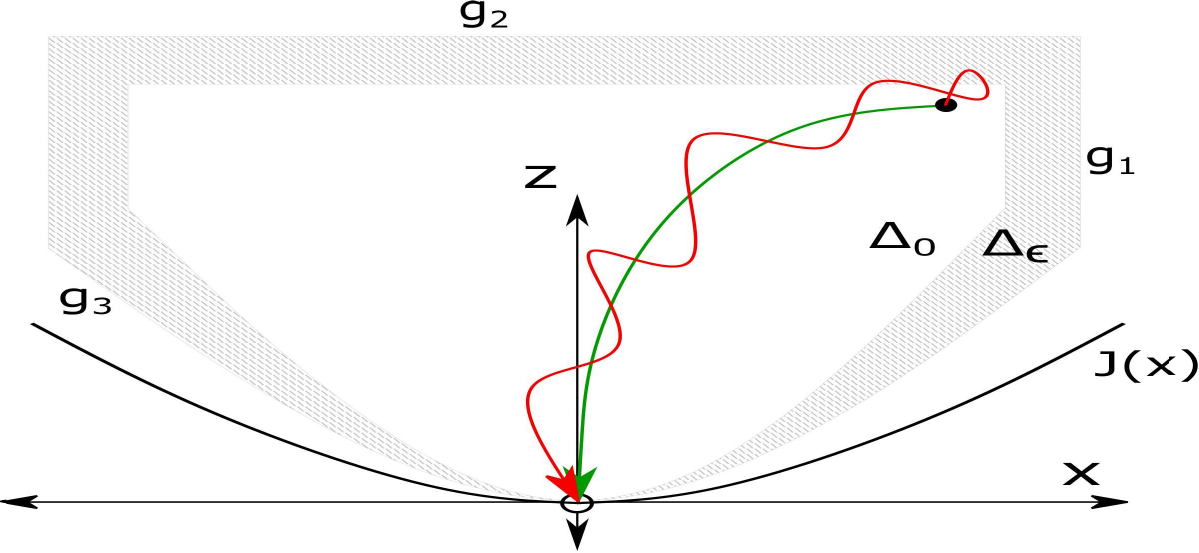}
        	\end{subfigure}
        	\begin{subfigure}[h]{0.49\columnwidth}
        	\centering
        	    \includegraphics[width=\textwidth,height=\textwidth]{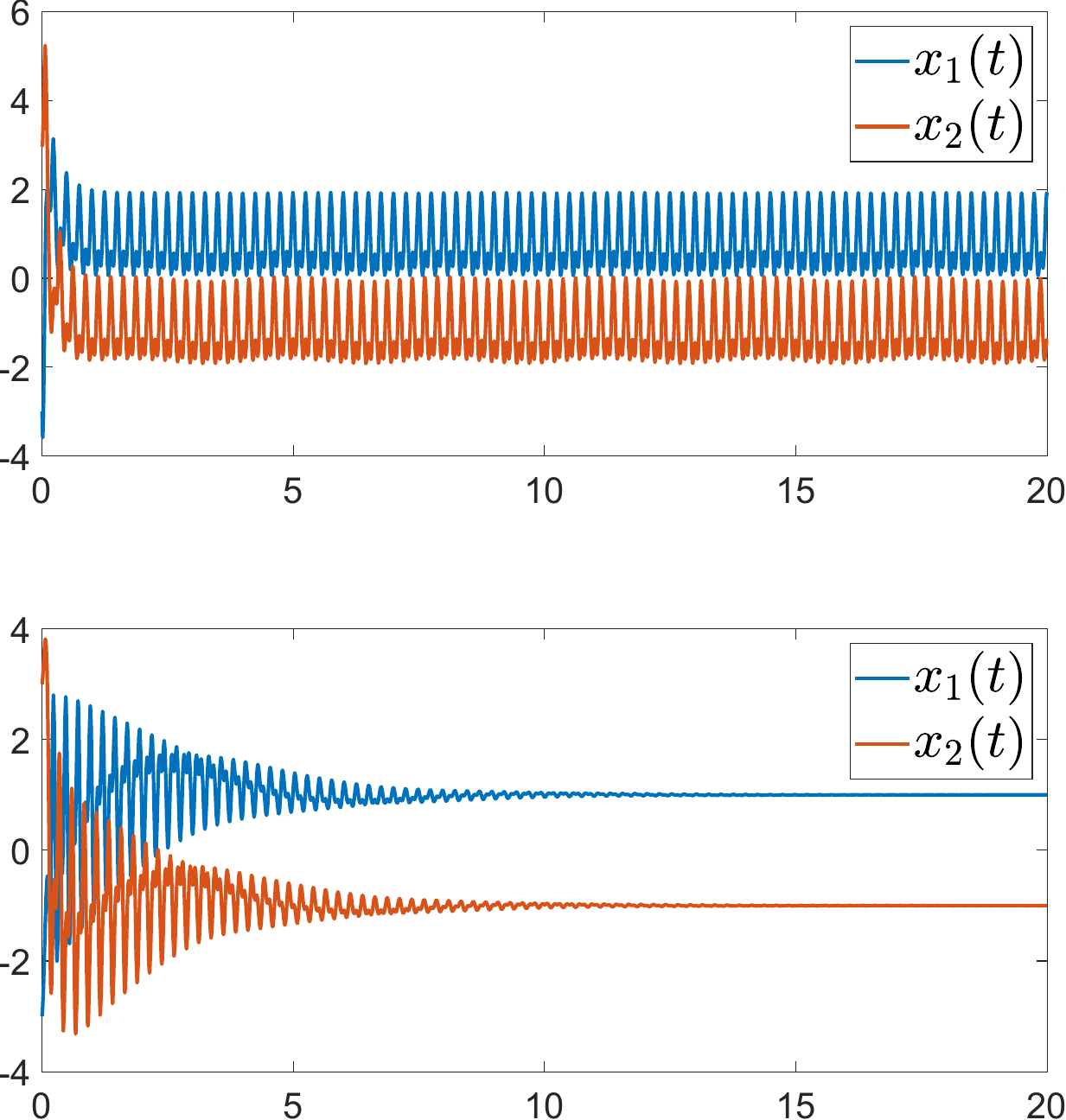}
        	\end{subfigure}
        	\caption{(left): Illustration of $\Delta_0,\Delta_\epsilon$ and sample trajectories, (right): Numerical results of Example 3.1 for the approach in \cite{Grushkovskaya2017} (top), and our approach (bottom)}
    	    \label{fig:exmp2}
	    \end{figure}
	\end{exmp}
	The numerical results for the proposed control law are shown in Fig.(\ref{fig:exmp2}), where we used the initial conditions $ \bm{x}(0) = (-3,3), ~z(0) = J(\bm{x}(0))+3 = 61,$ and the frequency parameters $\omega = 4,~ \omega_1 = 1,\omega_2 = 2$.
	
	\begin{rem}
	    We remark that the proposed method can tolerate bounded monotonic decrease of the minimum value of the function as demonstrated in the provided example. However, we emphasize that it does not tolerate general time-dependent variations of the cost function in the current formulation. This is due to the nature of the dynamic upper bound on the cost function (i.e. $z(t)$). 
	\end{rem}
\section{Conclusion and Future work}
 In this brief note, we propose an extension to extremum seeking control via Lie bracket approximations that allows asymptotic convergence to the minimum point for a cost function in the absence of information on its minimum value. The proposed control law leads to bounded control signals that vanish as the system converges to the minimum point, and bounded frequency of oscillation. We also a provide a procedure to obtain an estimate on the required frequency. Numerical simulations show that similar results may hold for the case of a dynamic cost function under appropriate assumptions on the dynamics.
	
	\section*{Acknowledgment}
	The authors like to acknowledge the support of the NSF Grant CMMI‐846308. The authors also thank the reviewers, whose suggestions helped to improve the manuscript. The first author thanks Prof. Anton Gorodetski for fruitful discussions, and Prof. Mostafa Abdallah for his continued support. 
	\begin{appendices}
	\section{Preliminary Results}
	Consider the Initial Value Problem (\textbf{IVP})
	\begin{equation}\label{IVP}
		\dot{\bm{\zeta}}(t) = \bm{f}_0(\bm{\zeta}(t)) + \sum\limits_{\lambda\in\Lambda} \bm{f}_{\lambda}(\bm{\zeta}(t)) ~u_{\lambda}(t),\quad \bm{\zeta}(0)\in \Xi_0
	\end{equation}
	where $\Xi_0\subset \Xi\subset\mathbb{R}^n$, $\Lambda$ is the set of all ordered pairs $(j,s), j\in\{1,2,...,n\},s\in\{1,2\}$, $\bm{f}_0,\bm{f}_\lambda\in \Gamma^2(\Xi)$ and the dither signals $u_\lambda(\cdot)$ are defined by Eq. (\ref{eq:dithers})
	\begin{lem}\thlabel{lem1}\cite{Durr2013,liu1997a,sutt2019}
	Let $g\in C^3(\Xi;\mathbb{R})$. Then, for every solution $\bm{\zeta}: I\rightarrow \Xi$ of (\ref{IVP}), the function $g\circ\bm{\zeta}: I\rightarrow \mathbb{R}$ satisfies
	\begin{equation*}
		\begin{aligned}
			g(\bm{\zeta}(t))\big|_{t_1}^{t_2} = R^g_1(\bm{\zeta}(t),t)\big|_{t_1}^{t_2} + \int\limits_{t_1}^{t_2} \big(F^g(\bm{\zeta}(t)) + R^g_2(\bm{\zeta}(t),t)\big)dt
		\end{aligned}
	\end{equation*}
	where $I$ is the interval of existence and uniqueness of $\bm{\zeta}(\cdot)$, $ t_1, t_2\in I, t_2>t_1$, and \\
	\begin{equation*}
	\begin{aligned}
	F^g(\bm{\zeta}) &= L_{\bm{f}_0}g(\bm{\zeta}) -  \sum\limits_{j=1}^m L_{[\bm{f}_{j,1},\bm{f}_{j,2}]}g(\bm{\zeta})\\
	R^g_1(\bm{\zeta},t) &=~ \sum\limits_{\lambda\in \Lambda}L_{\bm{f}_{\lambda}} g(\bm{\zeta})~U_{\lambda}(t) \\ &- \sum\limits_{\lambda_1,\lambda_2\in \Lambda}\hspace*{-0.075in}L_{\bm{f}_{\lambda_2}}L_{\bm{f}_{\lambda_1}} g(\bm{\zeta})~U_{\lambda_1,\lambda_2}(t)\\
	R^ g_2(\bm{\zeta},t) &= -\sum\limits_{\lambda\in \Lambda}L_{\bm{f}_0}L_{\bm{f}_{\lambda}} g(\bm{\zeta})~U_{\lambda}(t)\\ &+\sum\limits_{\lambda_1,\lambda_2\in \Lambda}L_{\bm{f}_{0}}L_{\bm{f}_{\lambda_2}}L_{\bm{f}_{\lambda_1}} g(\bm{\zeta})~U_{\lambda_1,\lambda_2}(t)\\&+\sum\limits_{\lambda_1,\lambda_2,\lambda_3\in \Lambda}\hspace*{-0.1in}L_{\bm{f}_{\lambda_3}}L_{\bm{f}_{\lambda_2}}L_{\bm{f}_{\lambda_1}} g(\bm{\zeta})~U_{\lambda_1,\lambda_2}(t) u_{\lambda_3}(t)
	\end{aligned}
	\end{equation*}
    \begin{equation*}
	\begin{aligned}
	U_{\lambda}(t) &= \int u_{\lambda}(t)~ dt \\
	U_{\lambda_1,\lambda_2}(t) &=\int \bigg(v_{\lambda_1,\lambda_2} + U_{\lambda_1}(t)~u_{\lambda_2}(t)\bigg) d\tau \\
	v_{\lambda_1,\lambda_2}&= \begin{cases}
	+1 & \lambda_1 = (j,1)\quad \& \quad \lambda_2 = (j,2) \\
	-1 & \lambda_1 = (j,2) \quad \& \quad \lambda_2 = (j,1) \\
	0  &\qquad\quad otherwise
	\end{cases}
	\end{aligned}
	\end{equation*}
	\end{lem}
	\begin{lem}\thlabel{lem2}
		Let $\Xi\subseteq\mathbb{R}^n, g_i\in C^3(\Xi;\mathbb{R})~\forall i\in\{1,2,...,r\}$. Let $\epsilon>0$ and define
		\begin{equation*}	
		\begin{aligned}
		\Delta_0 &= \big\{ \bm{\zeta}\in \Xi \big| g_i(\bm{\zeta})\leq 0, ~\forall i \in \{1,2,..~,r\} \big\}\\
		\Delta_\epsilon &= \big\{ \bm{\zeta}\in \Xi \big| g_i(\bm{\zeta})\leq\epsilon, ~\forall i \in \{1,2,..~,r\} \big\}
		\end{aligned}
		\end{equation*}
		and the subsets $ \Delta^i_\epsilon = \big\{ \bm{\zeta}\in \Delta_\epsilon \big| 0\leq g_i(\bm{\zeta})\leq\epsilon \big\}$. Suppose that $\forall i\in\{1,2,...,r\}$, whenever $\bm{\zeta} \in \Delta^i_\epsilon $, the following bounds hold
    	$$ \big\lVert R^{g_i}_1(\bm{\zeta},t)\big\lVert \leq \frac{c^{g_i}_1}{\sqrt{\omega}}, ~~~ \big\lVert R^{g_i}_2(\bm{\zeta},t)\big\lVert \leq \frac{c^{g_i}_2}{\sqrt{\omega}}, ~~~ F^{g_i}(\bm{\zeta}) \leq -b^{g_i}$$
	$\forall t\in\mathbb{R}$, where $c^{g_i}_1,c^{g_i}_2,b^{g_i}>0$ are constants. 
	Then $\exists\omega^*\in(0,\infty)$ such that $\forall\omega\in(\omega^*,\infty),\forall\bm{\zeta}(0)\in\Delta_0$ and maximal solution~ $\bm{\zeta}:\text{I}\rightarrow\Delta_\epsilon$ for the IVP (\ref{IVP}), where $0\in\text{I}=(t_e^-,t_e^+)$,$$\limsup\limits_{\tau\rightarrow t_e^+} g_i(\bm{\zeta}(\tau)) < \epsilon,~ \forall i\in\{1,2,..,r\}$$
	\end{lem}
	\begin{pf}
	Fix $\delta\in (0,\epsilon)$. If $g_i(\bm{\zeta}(t))<\delta, \forall t\in [0,t_e^+)$, the proof is complete. If not, then, by continuity of $g_i\circ\bm{\zeta}$ and the Intermediate Value Theorem, $\exists t_1,t_2\in I$, $t_2>t_1\geq0$, where $g_i(\bm{\zeta}(t_1)) = 0, ~ g_i(\bm{\zeta}(t_2))=\delta$, and $\bm{\zeta}(t)\in\Delta^i_\epsilon, ~\forall t\in [t_1,t_2]$.
	Using the bounds on $R^{g_i}_1,R^{g_i}_2, F^{g_i}$ and \thref{lem1}, we get
	$$g_i(\bm{\zeta}(t_2)) \leq \cancelto{0}{g_i(\bm{\zeta}(t_1))} \quad + \; \frac{2c^{g_i}_1}{\sqrt{\omega}} + \int\limits_{t_1}^{t_2}\bigg(-b^{g_i}+\frac{c^{g_i}_2}{\sqrt{\omega}}\bigg)~dt$$
	We define $$~\omega^* = \max\limits_{i\in\{1,2,...,r\}}\Big\{\Big(\frac{2 c_1^{g_i}}{\delta}\Big)^2, \Big(\frac{c_2^{g_i}}{b^{g_i}}\Big)^2\Big\}$$ and observe that $\forall \omega \in (\omega^*,\infty), \forall i\in\{1,2,...,r\}$, we have
	$$g_i(\bm{\zeta}(t_2)) < \delta \implies \limsup\limits_{\tau\rightarrow t_e^+} g_i(\bm{\zeta}(\tau)) < \epsilon \qed $$
	\end{pf}
	Note that $\omega^*$ in the proof of \thref{lem2} gives an estimate of the required frequency of oscillation in terms of the constants $c^{g_i}_1,c^{g_i}_2,b^{g_i}$ and a choice of $\delta\in(0,\epsilon)$. Thus, to choose a sufficiently large frequency, one needs to know the constants $c^{g_i}_1,c^{g_i}_2,b^{g_i}$. In the next lemma, we outline a procedure to estimate these constants under \thref{A1} for the static cost case. In fact, for a general dynamical system, if one can establish these bounds on the remainders, then the conclusions of \thref{main_thm} hold, as demonstrated by the numerical results provided above.
    
	\begin{lem}\thlabel{lem3}
	    Consider the system (\ref{es_sys}) and the functions (\ref{g_fcns}). Then, there exists constants $c_1^{g_i},c_2^{g_i} >0$ such that, $\forall \bm{\theta}\in\Delta^i_\epsilon, i \in\{1,2,3\}$ and $\forall t\in\mathbb{R}$:
	    $$ \big\lVert R^{g_i}_1(\bm{\theta},t)\big\lVert \leq \frac{c^{g_i}_1}{\sqrt{\omega}}, ~~~ \big\lVert R^{g_i}_2(\bm{\theta},t)\big\lVert \leq \frac{c^{g_i}_2}{\sqrt{\omega}}$$
	\end{lem}
	\begin{pf}
	    Via direct integration, the following bounds can be established 
	    \begin{equation*}
	        \begin{aligned}
	        &\left\lvert U_{\lambda_1} (t)\right\lvert\leq \frac{a_1}{\sqrt{\omega}},\quad  \left\lvert U_{\lambda_1,\lambda_2}(t)\right\lvert\leq \frac{a_1}{\omega}, \\
	        &\quad\quad \left\lvert U_{\lambda_1,\lambda_2}(t) u_{\lambda_3}(t)\right\lvert\leq \frac{a_1}{\sqrt{\omega}}
	        \end{aligned}
	    \end{equation*} 
	    $\forall \lambda_1,\lambda_2,\lambda_3\in\Lambda, \forall t\in\mathbb{R}$, where $a_1 > 0$ depends on the choice of the frequencies $\omega_j$. Due to space constraints, we only show how to establish a bound on one of the highest order terms in $ R^{g_3}_2(\cdot,\cdot)$. The rest of the bounds can be established following the same approach. Let $y = z-J(\bm{x})$. We compute
	{\begin{equation*}
	    \begin{aligned}
	        &L_{\bm{f}_{(i,1)}}L_{\bm{f}_{(i,1)}}L_{\bm{f}_{(i,1)}}g_3(\bm{\theta})= \partial^3_{i}g_3(\bm{\theta})F_1(y)^3 -2\partial^2_{i}g_3(\bm{\theta})\\ \times&F_1'(y)F_1(y)^2\partial_{i}J(\bm{x})
	        -\partial^2_{i}g_3(\bm{\theta})F_1'(y)F_1(y)^2\partial_i J(\bm{x}) \\ +&\partial_{i}g_3(\bm{\theta})F_1''(y)F_1(y)^2\partial_i J(\bm{x})^2 + \partial_{i}g_3(\bm{\theta})F_1'(y)^2  \\
	        \times&F_1(y)\partial_i J(\bm{x})^2-\partial_{i}g_3(\bm{\theta})F_1'(y)F_1(y)^2\partial^2_{i} J(\bm{x})
	    \end{aligned}
	\end{equation*}}
	It can be shown by direct computation that, for $y>0$, we have
	\begin{align*}
	    |F_1(y)| \leq &\sqrt{y}, ~ |F_1'(y)|\leq \frac{2}{\sqrt{y}}, ~ |F_2''(y)|\leq \frac{2}{y\sqrt{y}}\\
	    \partial_{i}g_3(\bm{\theta}) = &(2-\frac1m)\frac{\tilde{J}(\bm{x})^{1-\frac1m}\partial_iJ(\bm{x})}{z-J(x)}\text{sech}(\tilde{J}(\bm{x})^{2-\frac1m})^2\\
	    +&\frac{\tanh(\tilde{J}(\bm{x})^{2-\frac1m})\partial_iJ(\bm{x})}{(z-J(x))^2}\\
	    |\partial_i g_3(\bm{\theta})|\leq &\frac{2\tilde{J}(\bm{x})^{1-\frac1m}|\partial_iJ(\bm{x})|}{|z-J(x)|}+\frac{\tilde{J}(\bm{x})^{2-\frac1m}|\partial_iJ(\bm{x})|}{(z-J(x))^2}
	\end{align*}
	Using \thref{A1}, we can see that $$|\partial_iJ(\bm{x})|\leq\lVert \nabla J(\bm{x})\lVert\leq \sqrt{\kappa_2} \tilde{J}(\bm{x})^{1-\frac{1}{2m}}$$ Moreover, we know that $$\frac{\tanh(J(\bm{x})^{2-\frac1m})}{z-J(\bm{x})}\leq y_0+\epsilon,\quad \forall\bm{\theta}\in\Delta^3_\epsilon$$ Thus, it holds that
	\begin{equation*}
	    |\partial_i g_3(\bm{\theta})|\leq \sqrt{\kappa_2}(y_0+\epsilon)\frac{\tilde{J}(\bm{x})^{2-\frac{3}{2m}}}{\tanh(\tilde{J}(\bm{x})^{2-\frac1m})}(2+\tilde{J}(x))
	\end{equation*}
	For $\theta\in\Delta^3_\epsilon$, we have $\tilde{J}(\bm{x})\leq J_0$, and $|F'(y)F(y)|\leq 2$. Thus, we have:
	\begin{align*}
	    |\partial_i g_3(\bm{\theta})F'(y)F(y)\partial_i J(\bm{x})|\leq a_2 \frac{\tilde{J}(\bm{x})^{2-\frac{1}{m}}}{\tanh(\tilde{J}(\bm{x})^{2-\frac1m})}
	\end{align*}
	where $a_2=\kappa_2 J_0^{1-\frac{1}{m}}(y_0+\epsilon)(2+J_0)$. Finally, it can be shown that
	$$\frac{\tilde{J}(\bm{x})^{2-\frac{1}{m}}}{\tanh(\tilde{J}(\bm{x})^{2-\frac1m})}\leq\tilde{J}(\bm{x})^{2-\frac1m}+1$$
	This leads to the bound: $$ |\partial_i g_3(\bm{\theta})F'(y)F(y)\partial_i J(\bm{x})|\leq a_3 $$
	where $a_3 = a_2 (1+J_0^{2-\frac1m})$. Following a similar approach, it can be shown that all the terms in the Lie derivative are bounded,
	$|L_{\bm{f}_{(i,1)}}L_{\bm{f}_{(i,1)}}L_{\bm{f}_{(i,1)}}g_3(\bm{\theta})|\leq a_4$, where $a_4>0$ is the sum of all the bounds on the individual terms. Consequently, we have established the bound
	$$|L_{\bm{f}_{(i,1)}}L_{\bm{f}_{(i,1)}}L_{\bm{f}_{(i,1)}}g_3(\bm{\theta})U_{(i,1),(i,1)}(t) u_{(i,1)}(t)|\leq \frac{a_5}{\sqrt{\omega}}$$
	$\forall \bm{\theta}\in \Delta^3_\epsilon, \forall t\in\mathbb{R}$, where $a_5 = a_1a_4>0$. Following this procedure for each individual term in the remainders will give the explicit bounds on $R^{g_i}_1,R^{g_i}_2,i\in\{1,2,3\}$ in terms of the constants $\kappa_1,\kappa_2,\gamma_1,\gamma_2,\epsilon,J_0,y_0,z_0,\omega_j$.
	\end{pf}
	\section{Proof of Main Theorem}
	 \begin{pf}Let $J_0>0$ be such that the level set $$\mathcal{L}_{J_0} = \big\{\bm{x}\in\mathbb{R}^n\big|\tilde{J}(\bm{x})\leq J_0\big\} \subset D$$ Fix an $\epsilon> 0$, and let $$y_0 > \frac{1}{2\kappa_1}\big(1+\sqrt{1+8\kappa_1\epsilon}\big), ~~ z_0 > J_0+y_0$$ The functions $F_s(\cdot),s\in\{1,2\}$ are locally Lipschitz in $\Delta_\epsilon$. Hence, absolutely continuous maximal solutions of (\ref{es_sys}) with $\bm{\theta}(0)\in\Delta_\epsilon$ exist and are unique. We consider a maximal solution $\bm{\theta}: I \rightarrow \Delta_\epsilon$ of (\ref{es_sys}) with $\bm{\theta}(0)\in\Delta_0$ and apply \thref{lem1} to the functions $g_i,~i\in\{1,2,3\}$ defined by Eq. (\ref{g_fcns}). The next step is to establish the bounds on $F^{g_i}, R_1^{g_i}, R_2^{g_i}$ in \thref{lem2} for $g_i(\cdot), i\in\{1,2,3\}$. 
	We compute
	\begin{equation*}
	\begin{aligned}
	F^{g_1}(\bm{\theta}) &= -\left\lVert\nabla J(\bm{x})\right\lVert^2 \quad , \quad F^{g_2}(\bm{\theta}) = -z + J(\bm{x}) \\
	F^{g_3}(\bm{\theta}) &= \frac{\eta(\tilde{J}(\bm{x})^{2-\frac{1}{m}})}{(z-J(\bm{x}))}-\frac{\eta(\tilde{J}(\bm{x})^{2-\frac{1}{m}})\left\lVert\nabla J(\bm{x})\right\lVert^2}{(z-J(\bm{x}))^2}\\ &-\frac{\left(2-\frac{1}{m}\right) \tilde{J}(\bm{x})^{1-\frac{1}{m}}\eta'(\tilde{J}(\bm{x})^{2-\frac1m}) \left\lVert\nabla J(\bm{x})\right\lVert^2}{z-J(\bm{x})}
	\end{aligned}
	\end{equation*}
	where $\eta(y) = \tanh(y)$. We note that in case of $g_2$, the remainder terms $R^{g_2}_1, R^{g_2}_2$ in \thref{lem1} identically vanish, and the only remaining term inside the integral is $$F^{g_2}(\bm{\theta}) = -z + J(\bm{x}) < 0, ~ \forall \bm{\theta}\in \text{epi}_S(J)$$ We conclude, similar to the proof of \thref{lem2}, that  $g_2(\bm{\theta}(t)) \leq 0,~ \forall t\in I,~ \forall \omega \in (0,\infty)$.	Due to \thref{A1}, we know that $\forall\bm{\theta}\in\Delta^1_\epsilon$, we have $$F^{g_1}(\bm{\theta}) \leq - \kappa_1 J_0^{2-\frac{1}{m}}$$ Furthermore, by definition of $g_3(\cdot)$, and thanks to the property that 
	$\tanh(y) \leq y,~ \forall y\geq 0 $
	and the choice of $y_0$, we have $$F^{g_3}(\bm{\theta})\leq y_0+\epsilon-\kappa_1 y_0^2 < -\epsilon, \quad \forall\bm{\theta}\in\Delta^3_\epsilon$$ The bounds on the remainders $R^{g_2}_1, R^{g_2}_2$ can be explicitly computed, as outlined in \thref{lem3}. We now apply \thref{lem2} with the bounds established above to conclude that $\exists\omega^*\in(0,\infty)$ such that $\forall\omega\in(\omega^*,\infty),\forall\bm{\theta}(0)\in\Delta_0$ and maximal solution~ $\bm{\theta}:\text{I}\rightarrow\Delta_\epsilon$, where $0\in\text{I}=(t_e^-,t_e^+)$, we have $$\limsup\limits_{\tau\rightarrow t_e^+} g_i(\bm{\theta}(\tau)) < \epsilon,~\forall i\in\{1,2,3\}$$
	We note that the only remaining boundary in the definitions of $\Delta_0,\Delta_\epsilon$ is the point $(\bm{x}^*,J(\bm{x}^*))$. Clearly $\dot{z}(t)<0$. Moreover, we have that $$ \dot{z} = -z + J(\bm{x}) \geq -z + J(\bm{x}^*) \implies \tilde{z}(t)\geq \tilde{z}(0) e^{-t} >0$$ where $\tilde{z}(t) = z(t)-J(\bm{x}^*)$.
	Thus for any finite $t_e^+>0$, we have that $$\lim\limits_{\tau\rightarrow t_e^+} z(\tau)-J(\bm{x}(\tau)) > 0$$ This implies that $\forall\omega\in(\omega^*,\infty)$, maximal solutions that start inside $\Delta_0$ do not escape $\Delta_\epsilon$ in any finite time, hence $[0,\infty)\subset I$. Moreover, since z(t) is bounded below and strictly decreasing, we have that $$\lim\limits_{\tau\rightarrow +\infty} z(\tau)-J(\bm{x}(\tau)) = 0$$
	Consequently, we see that due to the definition of $g_3(\cdot)$, it must be true that $$\lim\limits_{\tau\rightarrow +\infty}\eta(\tilde{J}(\bm{x}(\tau))^{2-\frac{1}{m}}) = 0 \implies \lim \bm{x}(\tau)  = \bm{x}^*$$
	Combining all of the above, we conclude that $$\lim\limits_{\tau\rightarrow +\infty}\bm{\theta}(\tau) = (\bm{x}^*,J(\bm{x}^*))$$
	\end{pf}
	\end{appendices}
	\bibliographystyle{plain}
	\bibliography{main}
\end{document}